\input amstex
\nopagenumbers
\documentstyle{amsppt}
\refstyle{C}
%\magstep1
%\parindent=1em
%\baselineskip 15pt
%\hsize=12.3cm
%\vsize 18.5cm
\TagsOnRight
\define\e{@,@,@,@,@,}

\define\iv{\roman{Div}\e}
\define\br{\roman{Br}@,@,@,}
\define\g{\varGamma}
\define\ng{\negthickspace}
\define\cok{\roman{Coker}\,}
\define\ra{\rightarrow}
\define\krn{\roman{Ker}\,}
\define\img{{\roman{Im}}\,}
\define\be{@!@!@!@!@!}
\define\lbe{@!@!@!}
\define\kwstar{K\be\be @!@!\phantom{.}^{*}_{\ng w}}
\define\gal{\roman{Gal}\e}

\define\tprod{\tsize{\prod}}
\define\s{\Cal }

\font\tencyr=wncyr10
\font\sevencyr=wncyr7
\font\fivecyr=wncyr5
\newfam\cyrfam
\def\cyr{\fam\cyrfam\tencyr}
\textfont\cyrfam=\tencyr
\scriptfont\cyrfam=\sevencyr
\scriptscriptfont\cyrfam=\fivecyr
\document
\loadbold

\topmatter

\title Capitulation, ambiguous classes and the cohomology of the units.
\endtitle

\author Cristian D. Gonz\'alez-Avil\'es\endauthor

\address{Departamento de Matem\'aticas, Universidad
Andr\'es Bello, Rep\'ublica 252, Santiago, Chile.}
\endaddress

\email{cristiangonzalez\@unab.cl}\endemail

\bigskip

\rightheadtext{Capitulation, ambiguous classes and cohomology.}

\bigskip

\keywords{Ideal class groups, capitulation, units in global
fields, cohomology}\endkeywords

\thanks{The author is partially supported by Fondecyt research grant 1061209 and
Universidad Andr\'es Bello research grant DI-29-05/R.}\endthanks

\footnote"\phantom{.}"{2000 {\it Mathematics Subject
Classification:}\;\; 11R29, 11R27, 11R34.}

\abstract{This paper presents results on both the kernel and
cokernel of the $S$-capitulation map $C_{\be F,S}\ra C_{\be
K,S}^{\e G}$ for arbitrary finite Galois extensions $K/F$ of
global fields (with Galois group $G$) and arbitrary finite sets of
primes $S$ of $F$ (assumed to contain the archimedean primes in
the number field case).}

\endabstract

\endtopmatter

\heading 0. Introduction.\endheading

Let $K/F$ be a finite Galois extension of number fields with
Galois group $G$, and let $C_{F}$ and $C_{K}$ denote,
respectively, the ideal class groups of $F$ and $K$. The extension
of ideals from $F$ to $K$ induces a natural {\it capitulation map}
$j_{K/F}\: C_{F}\ra C_{K}^{\e G}$. An important problem in Number
Theory is to determine the kernel of $j_{K/F}$, which is usually
called the ``capitulation kernel". The classical Principal Ideal
Theorem of global class field theory (see [11, Theorem II.5.8.3,
p.168] for a generalized version of this theorem) asserts that
$\krn j_{K/F}$ is all of $C_{F}$ if $K$ is the Hilbert class field
of $F$. This fact motivated, quite early on, a great deal of
interest in the study of $\krn j_{\be K/\be F}$ for subfields $K$
of the Hilbert class field of $F$. As a result, most of the
existing literature on the Capitulation Problem is concerned with
the study of $\krn j_{K/F}$ for unramified abelian extensions $K$
of $F$ (see, e.g., [4, 12, 17, 25]) or, more generally, with the
kernel (and cokernel, in the case of [20]) of the
$S$-$T$-capitulation map $j_{\be K\be/\be F,S}^{\e T}\:
C_{F,S}^{\e T}\ra\left(C_{K,S}^{\e T}\right)^{\be G}$ for
$T$-tamely ramified and $S$-split abelian extensions $K/F$ (see
[19, 20] and [11, Corollary II.5.8.6, p.170]). Here $S$ denotes a
finite set of primes of $F$, which we always assume to contain all
archimedean primes in the number field case, and $T$ is a finite
set of non-archimedean primes of $F$ which is disjoint from $S$.
One of the few exceptions to this ``general rule" is the work of
B.Schmithals [23], who studied the kernel of $j_{K/F}$ for certain
types of possibly ramified cyclic extensions $K$ of a quadratic
number field $F$. However, the general problem of studying both
the kernel and cokernel of $j_{\be K\be/\be F,S}^{\e T}$ for {\it
arbitrary} finite Galois extensions $K$ of $F$ (i.e., not
necessarily $S$-split or with ramification locus equal to $T$) has
yet to be addressed. In particular, it is an unfortunate fact
that, despite the long history of research on the Capitulation
Problem, very little attention has been accorded the cokernel of
$j_{\be K\be/\be F,S}^{\e T}$. In effect, apart from the
contribution [20] already mentioned, the only other works known to
us which study the capitulation cokernel are [9, Appendix] and, in
the context of Iwasawa Theory, [13, 14, 15, 16, 18].

In this paper we study both the kernel and cokernel of the
$S$-capitulation map $j_{K/F,S}\: C_{F,S}\ra C_{K,S}^{\e G}$ for
arbitrary finite Galois extensions $K/F$ of global fields and
arbitrary sets $S$ as above. We show, for example, that $\krn
j_{K/F,S}$ may be identified with the subgroup of
$H^{1}(G,U_{K,S})$ of all cohomology classes which are locally
trivial at all places outside $S$ (See Corollary
2.5)\footnote{This result may be well-known to all specialists in
this area. See Remark 2.6.}. We also obtain a new generalization
of Hilbert's Theorem 94 (see Theorem 2.7), which applies to
possibly ramified cyclic extensions of global fields (it is an
open question to determine whether Theorem 2.7 below holds true
for arbitrary abelian extensions, i.e., whether there exists a
``ramified version" of Suzuki's Theorem [25]). In Section 3 we
obtain certain general results on $\cok j_{K/F,S}$ which have some
rather interesting consequences. For example, Theorem 3.7 states
that in the ``semisimple case" (i.e., when the degree of $K/F$ is
prime to the class number $h_{K,S}$) it is possible to determine
the structure of $H^{2}\be(G,U_{\be K,S}\be)$ {\it completely}.
Sections 4-8 give applications of the main results of the paper.
Theorem 4.1 states (roughly) that the structure of
$H^{1}(G,U_{K})$ is determined by that of $C_{F}$ and by the
ramification indices of $K/F$ if $F$ belongs to a certain class of
number fields and $K$ is equal to its own genus field relative to
$F$ (the proof of this result uses well-known theorems of
Tannaka-Terada, H.Furuya and C.Thi\'ebaud). Section 5 deals with
cyclic extensions and gives, under a certain hypothesis, a lower
bound for the number of ambiguous $S$-classes of $K$ which do not
come from $F$. See Theorem 5.2. This result may be regarded as a
``dual" of Hilbert's Theorem 94. The very brief Section 6 computes
the kernel and cokernel of the capitulation map for certain types
of imaginary extensions of function fields. An application to
imaginary Artin-Schreir extensions is given. Section 7 discusses
the case where $S$ is large relative to $K/F$, i.e., when $S$
contains all archimedean primes and all ramified primes of $K/F$
(we note that a significant portion of earlier work on the
Capitulation Problem has taken place in this setting). We show
that in this case the kernel of the capitulation map is naturally
isomorphic to $H^{1}(G,U_{K,S})$, and that its cokernel is a
certain group which measures the failure of the Hasse principle
for $H^{2}(G,U_{K,S})$. In Section 8, which concludes the paper,
we use results from [22] to show that the main theorems of
Sections 2 and 3 have natural analogs in the context of divisor
class groups. This Section also contains some results which
supplement those of [22]. Finally, the paper contains an Appendix
which relates certain integers that appear in the main text to the
ramification indices of $K/F$.

\smallpagebreak

Allow us to make here the following additional comments which may
help clarify the approach adopted in this paper. Let the global
field $F$ be given and let $S$ and $T$ be given finite sets of
primes of $F$ satisfying the conditions stated above. Then a natural
question is to study the $S$-$T$ capitulation map $j_{\be K\be/\be
F,S}^{\e T}$ for {\it varying} choices of $K$. If $K$ is chosen so
that all primes of $S$ split completely in $K$ and $K/F$ is (tamely)
ramified {\it exactly} over $T$, then we are in the setting of [20].
But other choices of $K$ are possible. For example, $K$ could be
chosen so that the set of primes of $F$ which ramify in $K$ is in
fact {\it disjoint} from $T$ (which is certainly the case in this
paper since here we consider $T=\emptyset$), but no conditions are
imposed regarding the behavior in $K$ of the primes in $S$. By
adopting this approach, we have been able to obtain results on the
cohomology of the units whenever information on Capitulation is
available (see, e.g., Theorem 4.1), and results on Capitulation
whenever information on the units is available (see, e.g., Theorem
5.2 and Example 5.3).

\smallpagebreak

One final comment is in order. Several of the results of this
paper immediately yield divisibility relations which involve the
class numbers $h_{F,S}$, $h_{K,S}$ and various other invariants
(see, e.g., Theorem 5.2). Since explicitly stating all such
divisibility relations would soon become quite repetitive, we have
decided not to state them at all. We are certain that the
interested readers will retrieve them without difficulty from the
corresponding statements (for an illustration, see Example 5.3).

\bigskip

\heading Acknowledgements.\endheading

I am deeply grateful to D.Lorenzini for sending me copies of
several of the references listed at the end of the paper, which
would otherwise have been unavailable to me. I also thank C.Maire
for sending me a copy of his thesis [19].

\heading 1. Settings and notations.\endheading

Let $F$ be a global field, i.e. $F$ is a finite extension of ${\Bbb
Q}$ (the ``number field case") or is finitely generated and of
transcendence degree 1 over a finite field of constants $k$ of
characteristic $p$ and cardinality $q$ (the ``function field case").
Let $K/F$ be a finite Galois extension of global fields with Galois
group $G$. In the function field case, we will write $k^{\e\prime}$
for the field of constants of $K$ and $q^{\e\prime}$ for its
cardinality. The infinite primes of a function field $F$ are the
primes which lie above the prime of $k(\e t\e)$ corresponding to the
pole of $t$. A function field having only one infinite prime is
called {\it{imaginary}}. Now let $S$ be any nonempty finite set of
primes of the global field $F$, containing the archimedean primes in
the number field case. Where confusion is unlikely, we will denote
by $S$ (also) the set of primes of $K$ which lie above the primes in
$S$. In the remaining instances, this set will be denoted by
$S_{K}$. The set of {\it non-archimedean} primes of $F$ which ramify
in $K$ will be denoted by $R$. Now, for each prime $v\in R\cup S$,
we fix once and for all a prime $w$ of $K$ lying above $v$. We will
write ${\s I}_{F,S}$ for the group of fractional ideals (or
divisors) of $F$ with support outside $S$. A similar notation will
apply to $K$. Finally, the set of archimedean primes of a number
field $F$ will be denoted by $S_{\infty}$.

\heading 2. The capitulation kernel.\endheading

We begin by considering the exact sequence of $G$-modules
$$
1\ra U_{K,S}\ra K^{*}\ra{\s I}_{K,S}\ra C_{K,S}\ra 0,
$$
which we split into two short exact sequences of $G$-modules as
follows:
$$
1\ra U_{K,S}\ra K^{*}\ra K^{*}\be/\e U_{K,S}\ra 1\tag 1
$$
and
$$
1\ra K^{*}\be/\e U_{K,S}\ra{\s I}_{K,S}\ra C_{K,S}\ra 0.\tag 2
$$

\medpagebreak

\proclaim{Lemma 2.1} We have
$$
H^{1}(G,{\s I}_{K,S})=0.
$$
\endproclaim
{\it Proof.} This is well-known. See, for example, [28, Lemma
2.1].\qed

\medpagebreak

Set $\br(K/F)=H^{2}(G,K^{*})=\krn\ng\left[\e\br F\ra(\br K)^{\e
G}\e\right]$. By the preceding lemma and Hilbert's Theorem 90, (1)
and (2) yield the following exact sequences:
$$
1\ra F^{*}\be/\e U_{F,S}\ra (K^{*}\be/\e U_{K,S})^{\e G}\ra
H^{1}(G,U_{K,S})\ra 0,\tag 3
$$
$$\align
0&\ra H^{1}(G,K^{*}\be/\e U_{K,S})\ra
H^{2}(G,U_{K,S})\ra\br(K/F)\\
&\ra H^{2}(G,K^{*}\be/\e U_{K,S})\ra H^{3}(G,U_{K,S}),\tag 4
\endalign
$$
$$\align
1&\ra (K^{*}\be/\e U_{K,S})^{\e G}\ra{\s I}_{K,S}^{\, G}\ra
(C_{K,S})^{\e G}_{\roman{trans}}\ra 0,\tag 5
\endalign
$$
$$
0\ra (C_{K,S})^{\e G}_{\roman{trans}}\ra C_{K,S}^{\e G}\ra
H^{1}(G,K^{*}\be/\e U_{K,S})\ra 0,\tag 6
$$
and
$$
0\ra H^{1}(G,C_{K,S})\ra H^{2}(G,K^{*}\be/\e U_{\be K,S})\ra
H^{2}(G,{\s I}_{\e K,S}),\tag 7
$$
where, by definition\footnote{The composite $C_{K,S}^{\e G}\!\ra\!
H^{1}(G,K^{*}\be/\e U_{K,S})\!\hookrightarrow\! H^{2}(G,U_{K,S})$ is
known as the {\it transgression map}, and $(C_{K,S})^{\e
G}_{\roman{trans}}$ might well be called the group of {\it
transgressive ambiguous classes}. However, $(C_{K,S})^{\e
G}_{\roman{trans}}$ is better known as the group of {\it strongly
ambiguous classes}.},
$$
(C_{K,S})^{\e G}_{\roman{trans}}=\krn\ng\left[\e C_{K,S}^{\e G}\ra
H^{1}(G,K^{*}\be/\e U_{K,S})\e\right].
$$
Note that, by (5), $(C_{K,S})^{\e G}_{\roman{trans}}$ is trivial
exactly when every ambiguous $S$-ideal of $K$ is principal.

For subsequent use, we note that the connecting homomorphism
$(K^{*}\be/\e U_{K,S})^{\e G}\ra H^{1}(G,U_{K,S})$ appearing in
(3) maps a class $\beta\e U_{K,S}\in (K^{*}\be/\e U_{K,S})^{\e G}$
to the cohomology class $\{\xi_{\e\beta}\}\in H^{1}(G,U_{K,S})$
which is represented by the 1-cocycle $\xi_{\e\beta}\:G\ra
U_{K,S}$ defined by $\xi_{\e
\beta}(\sigma)=\beta^{\e\sigma}\!/\beta\quad (\sigma\in G)$.

Now let $j_{\e K\!/\be F,S}\:C_{F,S}\ra C_{K,S}^{\e G}$ be the
natural capitulation map. It is not difficult to see, using the
general description of the connecting homomorphism $C_{K,S}^{\e
G}\ra H^{1}(G,K^{*}\be/\e U_{K,S})$ appearing in (6) (see, e.g.,
[2, p.97]), that the image of $j_{\e K\!/\be F,S}$ is contained in
$(C_{K,S})^{G}_{\roman{trans}}$. We will write $j^{\e\prime}_{\e
K\!/\be F,S}\:C_{F,S}\ra(C_{K,S})^{G}_{\roman{trans}}$ for the map
induced by $j_{\e K\!/\be F,S}$. Clearly, $\krn j^{\e\prime}_{\e
K\!/\be F,S}=\krn j_{\e K\!/\be F,S}$ and (6) immediately yields
the following proposition.

\proclaim{Proposition 2.2} There exists a natural exact sequence
$$
0\ra\cok j^{\e\prime}_{\e K\!/\be F,S}\ra\cok j_{\e K\!/\be
F,S}\ra H^{1}(G,K^{*}\be/\e U_{K,S})\ra 0.\qed
$$
\endproclaim

\medpagebreak

Now, we have a natural exact commutative diagram
$$
\minCDarrowwidth{.5cm} \CD 1 @>>> F^{*}\be/\e U_{F,S}
@>>>{\s I}_{F,S} @>>> C_{F,S} @>>>0\\
@.   @VVV     @VVV     @VV j^{\e\prime}_{\e K\!/\be
F,S}V @.\\
1@>>>(K^{*}\be/\e U_{\be K,S})^{G}@>>>{\s I}_{K,S}^{\e G}@>>>
(C_{K,S})^{G}_{\roman{trans}}@>>> 0\endCD\tag 8
$$
where the top row is the exact sequence (2) over $F$, the bottom row
is (5), and the left vertical map comes from (3). The middle
vertical map is injective and its cokernel has the following
well-known description.

\proclaim{Lemma 2.3} There exists a canonical isomorphism
$$
\cok\ng\left[\e{\s I}_{F,S}\ra\left({\s
I}_{K,S}\right)^{G}\e\right]\simeq\bigoplus_{v\in R\setminus S}
H^{1}\lbe(G_{w},U_{w}).
$$
\endproclaim
{\it Proof}. We give a proof of this well-known result because we
will need some of the maps defined below. Let $D=\sum_{v\notin
S}\sum_{w\lbe\mid\lbe v} n_{w}\be(D)\e w\in{\s I}_{K,S}$. Since $G$
permutes transitively the primes of $K$ lying above the same prime
of $F$, we have $D\in\left({\s I}_{K,S}\right)^{\be G}$ if and only
if, for every $v\notin S$, the coefficients $n_{w}\be(D)$ coincide
for all $w\ng\mid\ng v$. Write $n_{v}(D)$ for this common value.
Then $D$ belongs to the image of the map ${\s I}_{F,S}\ra {\s
I}_{K,S},\sum m_{v}v\mapsto\sum m_{v}\e e_{v}\be\sum_{w\lbe\mid\lbe
v}w$, if and only if $n_{v}(D)$ is divisible by $e_{v}$ for each
$v\notin S$. Thus there exists a canonical isomorphism
$$
\cok\ng\left[\e{\s I}_{F,S}\ra\left({\s
I}_{K,S}\right)^{G}\e\right]\simeq\bigoplus_{v\in R\setminus S}\Bbb
Z/e_{v}\e\Bbb Z
$$
which maps the class of $D\in\left({\s I}_{K,S}\right)^{\be G}$ to
$(\e n_{v}\be(\be D\be)\ng\mod e_{v})_{v\in R\setminus S}\in
\bigoplus_{v\in R\setminus S}\Bbb Z/e_{v}\e\Bbb Z$. On the other
hand, for each $w\notin S_{K}$, there exists a natural exact
sequence
$$
1\ra U_{w}\ra K_{w}^{*}\overset\roman{ord}_{w}\to\longrightarrow\Bbb
Z\ra 0,
$$
where the map $\roman{ord}_{w}$ assigns the value 1 to a fixed
uniformizing parameter of $K_{w}$. Let $v$ be the prime of $F$ lying
below $w$ and write $G_{w}=\gal(K_{w}/F_{v})$. We have a natural
exact commutative diagram
$$
\minCDarrowwidth{.6cm} \CD 1 @>>> U_{v}
@>>>F_{v}^{*}@>>>e_{v}\e\Bbb Z @>>>0\\
@.   @\vert     @\vert    @VV\subseteq V @.@.\\
1@>>>U_{\be w}^{\e G_{w}}@>>>(K_{w}^{*})^{G_{w}}@>>> \Bbb Z @>>>
H^{1}(G_{w},U_{w})@>>>0.\endCD
$$
Thus we have a canonical isomorphism
$$
\Bbb Z/e_{v}\e\Bbb Z\simeq H^{1}(G_{w},U_{w})
$$
which maps a class $\e m\ng\be\mod\! e_{v}\in\Bbb Z/e_{v}\e\Bbb Z$
to the cohomology class $\{\xi_{\e m}\}\in H^{1}(G_{w},U_{w})$
represented by the 1-cocycle $\xi_{\e m}\:G_{w}\ra U_{w}$ given by
$\xi_{\e m}(\sigma)=\beta^{\e\sigma}\!/\beta\quad(\sigma\in G_{w})$,
where $\beta\in K_{w}^{*}$ satisfies $\roman{ord}_{w}(\beta)=m$.\qed

\medpagebreak

We now apply the snake lemma to diagram (8) and obtain an exact
sequence
$$\align
0 &\ra\krn j_{\e K\!/\be F,S}\ra\cok\ng\left[\e F^{*}\be/\e
U_{F,S}\ra (K^{*}\be/\e U_{K,S})^{\e
G}\e\right]\overset{l}\to\longrightarrow \cok\ng\left[\e{\s
I}_{F,S}\ra\left({\s
I}_{K,S}\right)^{G}\e\right]\\
&\ra\cok j^{\e\prime}_{\e K\!/\be F,S}\ra 0,
\endalign
$$
where the map $l$ is induced by the natural map $K^{*}\be/\e
U_{K,S}\ra {\s I}_{K,S}$. This map fits into a commutative diagram
$$
\minCDarrowwidth{.6cm} \CD \cok\ng\left[\e F^{*}\be/\e U_{F,S}\ra
(K^{*}\be/\e U_{K,S})^{\e G}\e\right]@>l>> \cok\ng\left[\e{\s
I}_{F,S}\ra\left({\s I}_{K,S}\right)^{G}\e\right]\\
@VV\simeq V  @VV\simeq V\\
H^{1}\lbe(G,U_{K,S})@>\lambda>>\bigoplus_{v\in R\setminus S}
H^{1}\be(G_{w},U_{w}),
\endCD
$$
where the right vertical map is the isomorphism of Lemma 2.3, the
left vertical map is induced by the connecting homomorphism
$(K^{*}\be/\e U_{K,S})^{\e G}\ra H^{1}(G,U_{K,S})$ described
earlier, and the bottom (``localization") map $\lambda$ may be
described as follows: let $c\in H^{1}\be(G,U_{K,S})$ be
represented by the 1-cocycle $\xi\:G\ra
U_{K,S},\sigma\mapsto\beta^{\e\sigma}\!/\beta$, where $\beta\e
U_{K,S}\in (K^{*}\be/\e U_{K,S})^{\e G}$. Then the $v$-component
of $\lambda(c)$ ($v\in R\setminus S$) is the cohomology class in
$H^{1}\be(G_{w},U_{w})$ represented by the 1-cocycle $\xi_{\e
v}\:G_{w}\ra U_{w}$ given by $\xi_{\e
v}(\sigma)=\beta^{\e\sigma}\!/\beta\,\,(\sigma\in G_{w})$.

\smallpagebreak

The above argument yields the following result.

\proclaim{Theorem 2.4} There exists an exact sequence
$$
0\ra\krn j_{\e K\!/\be F,S}\ra H^{1}\lbe(G,U_{K,S})
\overset{\lambda}\to\longrightarrow \bigoplus_{v\in R\setminus S}
H^{1}\be(G_{w},U_{w})\ra\cok j^{\e\prime}_{\e K\!/\be F,S}\ra 0,
$$
where $j_{\e K\!/\be F,S}$ is the capitulation map,
$j^{\e\prime}_{\e K\!/\be F,S}\: C_{F}\ra (C_{K,S})^{\e
G}_{\roman{trans}}$ is induced by $j_{\e K\!/\be F,S}$ and
$\lambda$ is the localization map described above.
\endproclaim

By the description of the map $\lambda$ given above and the proof of
Lemma 2.3, the following is an immediate consequence of the theorem.

\proclaim{Corollary 2.5} The capitulation kernel $\krn j_{\e
K\be/\be F,S}$ is canonically isomorphic to the subgroup of
$H^{1}(G,U_{K,S})$ of all cohomology classes which are represented
by a 1-cocycle $\xi\:G\ra U_{K,S}$ of the form
$\xi(\sigma)=\beta^{\e\sigma}\be/\beta\,\,(\sigma\in G)$, where
$\beta\e U_{K,S}\in (K^{*}\be/\e U_{K,S})^{\e G}$ satisfies
$\roman{ord}_{w}(\beta)\equiv 0\pmod{e_{v}}$ for all $v\in
R\setminus S$.\qed
\endproclaim

\medpagebreak

{\bf{Remark 2.6.}} An equivalent form of the exact sequence of
Theorem 2.4 was previously obtained by H.van der Wall [28, proof of
Theorem 1.3, bottom of p.7] in the case that $F$ is a number field
and $S$ is the set of all archimedean primes of $F$. See also [23,
Theorem 2, p.46]. Further, we invite the reader to compare Corollary
2.5 (for number fields and $S=S_{\infty}$) with [23, Corollary,
p.46].

\medpagebreak

Now set
$$
S^{\e\prime}=S\cup (R\setminus S\e)=S\cup R.
$$
We define integers $d_{v}$, for $v\in S^{\e\prime}$, when $K/F$ is a
{\it cyclic} extension of degree $n$ as follows:
$$
d_{v}=\cases [K_{w}\:\be F_{v}]\quad\text{ if $v\in S$}\\
\,\,\,\,\,\,\,\,\,e_{v}\qquad\,\,\e\text{ if $v\in R\setminus S$}.
\endcases\tag 9
$$
Clearly, each $d_{v}$ is a divisor of $n$. Set

$$\align
D&=\roman{l.c.m.}\{d_{v}\:v\in S^{\e\prime}\}\tag 10\\
\intertext{and} n_{0}&=n/\e D.\tag 11
\endalign
$$
Then we have the following generalization of Hilbert's Theorem 94.

\proclaim{Theorem 2.7} Let $F$ be a global field and let $K/F$ be
a cyclic extension of degree $n$. Let $d_{v}$, $D$ and $n_{0}$ be
the integers (9), (10) and (11), respectively. Then at least
$n_{0}\e\big/\lbe\be\left(n_{0},\tprod_{\e v\in
S^{\e\prime}}d_{v}/ D\e\right)$ $S$-ideal classes of $F$
capitulate in $K$.
\endproclaim
{\it Proof.} The proof of Lemma 2.3 and Theorem 2.4 immediately
yield the order relation
$$
\left[\krn j_{\e K\!/\be F,S}\right]\e\tprod_{\e v\in R\setminus
S}\e e_{v}=\left[H^{1}\lbe(G,U_{K,S})\right]\e\big[\cok
j^{\e\prime}_{\e K\!/\be F,S}\e\big].
$$
On the other hand, the well-known formula for the Herbrand
quotient of the $G$-module $U_{K,S}$ (see, e.g., [26, Proof of
Theorem 8.3, p.178]) yields the identity
$$
n \big[\widehat{H}^{\e
0}\lbe(G,U_{K,S})\big]=\left[H^{1}\lbe(G,U_{K,S})\right]\e
\tprod_{\e v\in S}\e[K_{w}\:\be F_{v}].
$$
Combining the preceding formulas, we obtain the equality
$$
\frac{n}{\left(n,\tprod_{\e v\in S^{\e\prime}}d_{v}\right)}
\big[\widehat{H}^{\e 0}\lbe(G,U_{K,S})\big]\e\big[\cok
j^{\e\prime}_{\e K\!/\be F,S}\e\big]=\frac{\tprod_{\e v\in
S^{\e\prime}}\e d_{v}}{\left(n,\tprod_{\e v\in
S^{\e\prime}}d_{v}\right)}\left[\krn j_{\e K\!/\be F,S}\right].
$$
It follows that $n\big/\be\left(n,\tprod_{\e v\in
S^{\e\prime}}d_{v}\right)=n_{0}\big/\be\left(n_{0},\tprod_{\e v\in
S^{\e\prime}}d_{v}/ D\e\right)$ divides $\left[\krn j_{\e K\!/\be
F,S}\right]$, as asserted.\qed

\medpagebreak

The group $\cok j^{\e\prime}_{\e K\!/\be F,S}$ appears to be as
difficult to compute as $\krn j_{\e K\!/\be F,S}\,$\footnote{ See
[28] for the computation of $\big[\e\cok j^{\e\prime}_{\e K\!/\be
F,S}\e\big]$ when $S=S_{\infty}$ in some particular cases, notably
when $K/F$ is an extension of prime degree of a quadratic number
field $F$.}. We close this section by giving an alternative
description of $\cok j^{\e\prime}_{\e K\!/\be F,S}$ (see Proposition
2.9 below) which may prove useful in future research on this group.
We need the following approximation lemma.

\proclaim{Lemma 2.8} The natural map
$$
F^{*}\be/\e U_{F,S}\ra\bigoplus_{v\in R\setminus S}F_{v}^{*}/\e
U_{v}
$$
is surjective.
\endproclaim
{\it Proof.} For each $v\in R\setminus S$, let $x_{v}U_{v}\in
F_{v}^{*}/\e U_{v}$ and set
$m=\roman{max}\{\roman{ord}_{v}(x_{v})\:v\in R\setminus S\}$. By
the strong approximation theorem [2, p.67], there exists a
$\beta\in F^{*}$ such that $\roman{ord}_{v}(\beta-x_{v})>m$ for
all $v\in R\setminus S$. It follows that
$\roman{ord}_{v}(\beta)=\roman{ord}_{v}(x_{v})$ for every $v\in
R\setminus S$, i.e., $\beta\e x_{v}^{-1}\in U_{v}$ for every $v\in
R\setminus S$. Consequently $\beta\e U_{F,S}\in F^{*}\be/\e
U_{F,S}$ maps to $(x_{v}U_{v})_{v\in R\setminus
S}\in\bigoplus_{v\in R\setminus S}F_{v}^{*}/\e U_{v}$.\qed

\medpagebreak

We now consider the exact commutative diagram
$$
\minCDarrowwidth{.5cm} \CD 1 @>>> F^{*}\be/\e U_{F,S}
@>>>(K^{*}\be/\e U_{K,S})^{\e G} @>>> H^{1}(G,U_{K,S})@>>>0\\
@.   @VVV     @VVV     @VV \lambda V @.\\
1@>>>\underset{v\in R\setminus S}\to{\bigoplus}F_{v}^{*}/\e
U_{v}@>>>\underset{v\in R\setminus S}\to{\bigoplus}(K_{w}^{*}/\e
U_{w})^{\e G_{w}}@>>>\underset{v\in R\setminus S}\to{\bigoplus}
H^{1}\be(G_{w},U_{w})@>>> 0,\endCD
$$
where the top row is the exact sequence (3), the bottom row is the
direct sum over $v\in R\setminus S$ of exact sequences analogous
to (3) for the extension $K_{w}/F_{v}$, and the unlabeled vertical
maps are the natural ones. By the preceding lemma, the above
diagram yields an isomorphism
$$
\cok\lambda\simeq\cok\ng\left[(K^{*}\be/\e U_{K,S})^{\e G}\ra
\underset{v\in R\setminus S}\to{\bigoplus}(K_{w}^{*}/\e U_{w})^{\e
G_{w}}\right].\tag 12
$$
Let
$$
\varphi\:(K^{*}\be/\e U_{K,S})^{\e G}\ra\bigoplus_{v\in R\setminus
S}\Bbb Z
$$
be given by
$$
\varphi(x\e U_{K,S})=(\roman{ord}_{w}(x))_{\e v\e\in\e R\setminus
S} \qquad(x\in K^{*}), \tag 13
$$
i.e., $\varphi$ is the composite of the natural map $(K^{*}\be/\e
U_{K,S})^{\e G}\ra\bigoplus_{v\in R\setminus S}(K_{w}^{*}/\e
U_{w})^{\e G_{w}}$ and the isomorphism
$\bigoplus\roman{ord}_{w}\:\bigoplus_{v\in R\setminus
S}(K_{w}^{*}/\e U_{w})^{\e
G_{w}}\overset{\sim}\to{\longrightarrow}\bigoplus_{v\in R\setminus
S}\Bbb Z$. Then (12) induces an isomorphism
$$
\cok\lambda\simeq\cok\varphi.
$$
On the other hand, Theorem 2.4 yields an isomorphism $\cok
j^{\e\prime}_{\e K\!/\be F,S}\simeq\cok\lambda$. Thus the
following holds.

\proclaim{Proposition 2.9} There exists an isomorphism
$$
\cok j^{\e\prime}_{\e K\!/\be F,S}\simeq\cok\varphi,
$$
where $\varphi$ is the map (13).\qed
\endproclaim

\medpagebreak

\heading 3. The capitulation cokernel.\endheading

We now recall the exact sequence (4):
$$\align
0&\ra H^{1}(G,K^{*}\be/\e U_{\be K,S})\ra H^{2}(G,U_{\be
K,S})\ra\br(K/F)\ra H^{2}(G,K^{*}\be/\e
U_{\be K,S})\\
&\ra H^{3}(G,U_{\be K,S}).
\endalign
$$
We note that the map $H^{2}(G,U_{\be K,S}\be)\ra\br(K/F)=H^{@,@,
2}(G,K^{*})$ appearing above is induced by the inclusion $U_{\be
K,S}\subset K^{*}$. Now let
$$
\roman{B}@,@,@,@,({\Cal O}_{F,S},{\Cal O}_{K,S})=
\krn\ng\left[\,\br(K/F)\ra H^{2}(G,{\Cal I}_{K,S})\e\right],\tag
14
$$
where the map involved is induced by the natural map $K^{*}\ra{\Cal
I}_{K,S},\, x\ra(x)$. Then, by (7), the above exact sequence induces
an exact sequence\footnote{ This exact sequence was discovered
independently by M.Auslander and A.Brumer in [3] and by S.Chase,
D.Harrison and A.Rosenberg in [6]. It was used in [5] and in [16].
More recently, it resurfaced in [9, Appendix], where a particular
case of (15) was derived from a Hochschild-Serre spectral sequence
in \'etale cohomology.}
$$\align
0&\ra H^{1}(G,K^{*}\be/\e U_{\be K,S}\be)\ra H^{2}(G,U_{\be
K,S}\be)\ra\roman{B}@,@,@,@,({\Cal O}_{\be
F,S},{\Cal O}_{\be K,S})\ra H^{1}(G,C_{K,S})\\
&\ra H^{3}(G,U_{\be K,S}\be).\tag 15
\endalign
$$
See [16, \S 1] for the details. By combining the preceding exact
sequence with Proposition 2.2, we obtain the following result.

\proclaim{Proposition 3.1} There exists a natural six-term exact
sequence
$$\align
0&\ra \cok j^{\e\prime}_{\e K\!/\be F,S}\ra\cok j_{\e K\!/\be
F,S}\ra H^{2}(G,U_{\be K,S}\be)\ra\roman{B}({\Cal O}_{\be
F,S},{\Cal O}_{\be K,S})\\
&\ra H^{1}(G,C_{K,S})\ra H^{3}(G,U_{\be K,S}\be),
\endalign
$$
where $\roman{B}({\Cal O}_{\be F,S},{\Cal O}_{\be K,S})$ is the
group (14).\qed
\endproclaim

\bigpagebreak

We will now use class field theory to compute $\roman{B}({\Cal
O}_{\be F,S},{\Cal O}_{\be K,S})$. This computation generalizes
[9, Lemma A.1].

\smallpagebreak

By global class field theory, $\br\e F$ is naturally isomorphic to
the kernel of the map $\sum\roman{inv}_{v}\:\bigoplus_{\text{all
$\be v$}}\br\e F_{\be v}\ra\Bbb Q/\e\Bbb Z$. Under this
isomorphism, the subgroup $\br(K/F)$ of $\br\e F$ may be
identified with the kernel of the induced map
$$
\sum\e\roman{inv}_{v}\:\bigoplus_{\text{all $\be v$}}
H^{2}(G_{w},\kwstar)\ra\Bbb Q/\e\Bbb Z,
$$
where, for each $v$, we regard $H^{2}\be(G_{w},K_{w}^{*})$ as a
subgroup of $\br\e F_{\be v}$ via the inflation map
$\roman{Inf}_{w/v}\: H^{2}\be(G_{w},K_{w}^{*})\ra\br(F_{v})$
(which is injective). On the other hand, there exist well-known
canonical isomorphisms
$$\align
H^{2}(G,{\Cal I}_{K,S})&\simeq\bigoplus_{v\notin
S}H^{2}\be(G_{w},\langle\,{\frak p}_{w}\rangle)
\simeq\bigoplus_{v\notin S}H^{2}\be(G_{w},\Bbb Z)\\
&\simeq \bigoplus_{v\notin S}H^{2}\be(G_{w},\kwstar/\e U_{w}).
\endalign
$$
Let $\psi\:H^{2}(G,{\Cal I}_{K,S})\ra\bigoplus_{v\notin
S}H^{2}\be(G_{w},\kwstar/\e U_{w})$ be the composite of the above
isomorphisms. Then the diagram
$$
\minCDarrowwidth{.6cm} \CD \br(K/F) @>>>\left(\e\underset{v\in
S}\to\bigoplus
H^{2}(G_{w},\kwstar)\right)\!\oplus\!\left(\underset{v\notin
S}\to\bigoplus H^{2}(G_{w},\kwstar)\right)
\\
@VVV  @VV\left(f_{1},\,\bigoplus_{v\notin S}\! f_{\e 2,v}\right)V\\
H^{2}(G,{\Cal I}_{K,S})@>\psi>>\underset{v\notin S}\to\bigoplus\,
H^{2}\lbe(G_{w},\kwstar/\e U_{w}),
\endCD
$$
where $f_{1}$ is the zero map and
$f_{2,v}\:H^{2}(G_{w},\kwstar)\ra H^{2}(G_{w},\kwstar/\e U_{w})$,
for each $v\notin S$, is induced by the natural map
$\kwstar\ra\kwstar/\e U_{w}\e$, commutes. It follows that the
identification
$$
\br(K/F)=\krn\ng\left[\,\,\bigoplus_{\text{all } v}
H^{2}\lbe\be\left(G_{w},\kwstar\right)\overset{\sum\roman{inv}_{\be
v}}\to{\longrightarrow} \Bbb Q/\e\Bbb Z\e\right]
$$
induces an identification
$$
\roman{B}@,@,@,({\Cal O}_{\be F,S},{\Cal O}_{\be
K,S})=\krn\ng\left[\,\,\bigoplus_{v\in S}\be
H^{2}(G_{w},\kwstar)\oplus\bigoplus_{v\notin S}\krn f_{\e 2,v}
\overset{\sum\roman{inv}_{\be v}}\to{\longrightarrow}\Bbb Q/\e\Bbb
Z\e\right].\tag 16
$$
Now, since $H^{1}(G_{w},\kwstar/\e U_{w})\simeq H^{1}(G_{w},\Bbb
Z)=0$ for each $v\notin S$, there exist canonical isomorphisms
$$
\krn f_{\e 2,v}\simeq H^{2}(G_{w},U_{w})\quad(v\notin S\e). \tag
17
$$
Note that the latter group is zero if $v\notin R$. Now recall the
set $S^{\e\prime}=S\cup (R\setminus S)=S\cup R$. We define integers
$d_{v}$, for $v\in S^{\e\prime}$, as follows\footnote{If $K/F$ is
cyclic, then the integers (18)-(19) agree with the integers
(9)-(10). See, for example, [1].}:
$$
d_{v}=\cases\,\,\,\,\,\,[K_{w}\:\be F_{v}]\,\,\,\,\qquad\text{ if
$v\in S$}\\
[\e H^{2}(G_{w},U_{w})\e]\quad\,\,\e\text{ if $v\in R\setminus
S$}.
\endcases\tag 18
$$
Further, set
$$
D=\roman{l.c.m.}\{d_{v}\:v\in S^{\e\prime}\}.\tag 19
$$
Then the invariant map $\roman{inv}_{v}$ induces isomorphisms
$$\align
H^{2}(G_{w},\kwstar)&\simeq d_{v}^{\e -1}\Bbb Z/\Bbb Z \qquad(v\in
S\e)\tag 20\\
H^{2}(G_{w},U_{w})&\simeq d_{v}^{\e -1}\Bbb Z/\e\Bbb Z\qquad(v\in
R\setminus S).\tag 21
\endalign
$$
It follows from (16), (17), (20) and (21) that there exists a
natural isomorphism
$$
\roman{B}@,@,@,({\Cal O}_{\be F,S},{\Cal O}_{\be
K,S})\simeq\krn\ng\be\left[\,\,\bigoplus_{v\in S^{\e\prime}}\e
d_{v}^{\e -1}\Bbb Z/\e\Bbb Z\overset{\Sigma}\to{\longrightarrow}
D^{-1}\Bbb Z/\Bbb Z\e\right],
$$
where $\Sigma$ is the summation map $(x_{v})\ra\sum\e x_{v}$. The
latter map is surjective (see [10, Lemma 1.2]), whence the
following holds.

\proclaim{Lemma 3.2} There exists an exact sequence
$$
0\ra\roman{B}@,@,@,({\Cal O}_{\be F,S},{\Cal O}_{\be K,S})\ra
\bigoplus_{v\in S^{\e\prime}}\e\Bbb Z/\e d_{v}\e\Bbb
Z\overset{\Sigma}\to{\longrightarrow}\Bbb Z/ D\e\Bbb Z\ra 0,
$$
where $d_{v}$ and $D$ are the integers (18) and (19),
respectively. In particular,
$$
\left[\,\roman{B}@,@,@,({\Cal O}_{\be F,S},{\Cal O}_{\be K,S})
\,\right]=\left(\e\tprod_{\e v\in
S^{\e\prime}}d_{v}\right)\lbe/D.\qed
$$
\endproclaim

\medpagebreak

Now, combining Proposition 3.1 and Lemma 3.2, we obtain the
following result.

\proclaim{Theorem 3.3} There exists an exact sequence
$$\align
0&\ra \cok j^{\e\prime}_{\e K\!/\be F,S}\ra\cok j_{\e K\!/\be
F,S}\ra H^{2}(G,U_{\be K,S}\be)\ra\roman{B}({\Cal O}_{\be
F,S},{\Cal O}_{\be K,S})\\
&\ra H^{1}(G,C_{K,S})\ra H^{3}(G,U_{\be K,S}\be),
\endalign
$$
where $\roman{B}({\Cal O}_{\be F,S},{\Cal O}_{\be K,S})$ is a
group of order $\left(\e\tprod_{\e v\in
S^{\e\prime}}d_{v}\right)\lbe/D$.\qed
\endproclaim

\medpagebreak

{\bf Remark 3.4.} We note that, via the isomorphisms (16) and
(17), the map $H^{2}(G,U_{K,S})\ra \roman{B}@,@,@,({\Cal O}_{\be
F,S},{\Cal O}_{\be K,S})$ appearing in the exact sequence of the
theorem is induced by the inclusions $U_{K,S}\subset U_{w}$
($v\notin S$) and $U_{K,S}\subset\kwstar$ ($v\in S$). Thus, by
(15),
$$\align
H^{1}(G,K^{*}\be/\e U_{\be K,S}\be)&\simeq\krn\ng\left[\e
H^{2}(G,U_{K,S})\ra \roman{B}@,@,@,({\Cal O}_{\be F,S},{\Cal
O}_{\be K,S})\e\right]\\
&\simeq\krn\ng\be\left[H^{2}\be(G,U_{K,S})\be\ra\!\bigoplus_{v\in
R\setminus S}\! H^{2}(G_{w},U_{w})\lbe\oplus\lbe\bigoplus_{v\in
S}\be H^{2}(G_{w},\kwstar)\right],
\endalign
$$
where the last map is induced by the inclusions $U_{K,S}\subset
U_{w}$ ($v\notin S$) and $U_{K,S}\subset\kwstar$ ($v\in S$).

\medpagebreak

We now derive some consequences of Theorem 3.3. The first one
generalizes [22, Theorem 2(C), p.161] (see Proposition 2.2).

\medpagebreak

\proclaim{Corollary 3.5} Assume that the integers $d_{v}$, where
$v\in S^{\e\prime}$, are pairwise coprime. Then there exists a
natural exact sequence
$$
0\ra \cok j^{\e\prime}_{\e K\!/\be F,S}\ra\cok j_{\e K\!/\be
F,S}\ra H^{2}(G,U_{\be K,S}\be)\ra 0.\qed
$$
\endproclaim

\medpagebreak

\proclaim{Corollary 3.6} Assume that $H^{2}\be(G,U_{\be K,S}\be)=0$.
Then $H^{1}\be(G,C_{K,S})$ contains at least $\left(\e\tprod_{\e
v\in S^{\e\prime}}d_{v}\e\right)\lbe\big/D$ elements.\qed
\endproclaim

\bigpagebreak

The next theorem refers to the ``limit case" where $[\e K\be:F\e]$
is prime to $h_{K,S}\e$\footnote{This condition can always be
satisfied by enlarging $S$ appropriately. We also note that this
case has been previously studied by H.Yokoi in [29, \S 4].}.

\proclaim{Theorem 3.7} Let $K/F$ be a Galois extension of global
fields, of degree $n$. Assume that $n$ is prime to $h_{K,S}$. Then
there exists a natural isomorphism
$$
H^{2}\be(G,U_{\be K,S}\be)\simeq \krn\ng\left[\e \bigoplus_{v\in
S^{\e\prime}}\e\Bbb Z/\e d_{v}\e\Bbb
Z\overset{\Sigma}\to{\longrightarrow}\Bbb Z/ D\e\Bbb Z\e\right].
$$
In particular,
$$
\left[\e H^{2}\be(G,U_{\be K,S}\be)\e\right]=\left(\e\tprod_{\e
v\in S^{\e\prime}}d_{v}\right)\lbe/D.
$$
\endproclaim
{\it Proof.} Let $N\:C_{K,S}\ra C_{K,S}^{\e G}$ be the map
$c\mapsto\tprod_{\e\sigma\in G}\e c^{\e\sigma}$. Then $N=j_{\e
K\!/\be F,S}\circ N_{K/F}$, where $N_{K/F}\: C_{K,S}\ra C_{F,S}$
is induced by the norm of ideals. Consequently, $\cok j_{\e
K\!/\be F,S}$ is a quotient of $C_{K,S}^{\e
G}/NC_{K,S}=\widehat{H}^{0}(G,C_{K,S})$. Since the latter group is
annihilated by multiplication by $(n,h_{K,S})=1$, we conclude that
$\cok j_{\e K\!/\be F,S}=0$. The theorem now follows from Theorem
3.3 since $H^{1}(G,C_{K,S})$ is also annihilated by multiplication
by $(n,h_{K,S})=1$.\qed

\medpagebreak

{\bf Remark 3.8.} If $K/F$ is {\it cyclic} of degree $n$, and $n$
is prime to $h_{K,S}$, then the order of $\krn j_{\e K\!/\be F,S}$
may be computed explicitly. Indeed, by the theorem,
$\big[\e\widehat{ H}^{0}\be(G,U_{\be K,S}\be)\e\big]=\left[\e
H^{2}\be(G,U_{\be K,S}\be)\e\right]=\left(\e\tprod_{\e v\in
S^{\e\prime}}d_{v}\right)\lbe/D$. On the other hand, the proof of
the above theorem and Proposition 2.2 show that $\cok
j^{\e\prime}_{\e K\!/\be F,S}=0$. Finally the proof of Theorem 2.7
yields
$$
\left[\e\krn j_{\e K\!/\be F,S}\right]=n_{@,@,0},
$$
where $n_{@,@,0}=n/D$.

\heading 4. Genus fields.\endheading

In this section we consider abelian extensions of certain base
fields $F$ which are equal to their own genus field relative to $F$.
For such extensions (and a minimal set $S$), the group
$(C_{K,S})^{\e G}_{\roman{trans}}$ is zero and the following holds.

\proclaim{Theorem 4.1} Let $K/F$ be a finite abelian extension of
number fields. Assume that $K$ is its own genus field relative to
$F$. Assume, furthermore, that one of the following conditions
holds:

\roster

\item"(a)" $K/F$ is cyclic, or

\item"(b)" $F$ is either the rational field $\Bbb Q$ or an
imaginary quadratic extension of $\Bbb Q$ of discriminant $<-4$
and conductor prime to 2.
\endroster
Then there exists an exact sequence
$$
0\ra C_{F}\ra
H^{1}\lbe(G,U_{K})\ra\bigoplus_{v\,\text{non-arch.}}\ng\Bbb
Z/e_{v}\e\Bbb Z\ra 0.
$$
Furthermore, there exists a canonical isomorphism
$$
C_{K}^{\e G}\simeq H^{1}(G,K^{*}\be/\e U_{\be K,S}\be).
$$
In particular, if $h_{F}=1$, there exists an isomorphism\footnote{
More generally, if $C_{F}=C_{F}^{\e e_{\lbe v}}$ for every
non-archimedean prime $v$, then $H^{1}\lbe(G,U_{K})\simeq\!
C_{F}\oplus\left(\bigoplus_{v\text{ non-arch.}}\be\Bbb Z/e_{v}\e\Bbb
Z\right)$.}
$$
H^{1}\lbe(G,U_{K})\simeq\!\bigoplus_{v\text{ non-arch.}}\be\Bbb
Z/e_{v}\e\Bbb Z.
$$
\endproclaim
{\it Proof.} We apply results from the preceding sections with $S$
equal to the set of all archimedean primes of $F$. By theorems of
Tannaka-Terada, H.Furuya [8] and C.Thiebaud [27], every ambiguous
ideal of $K$ is principal. Therefore $(C_{K})^{\e
G}_{\roman{trans}}=0$ (see (5)), whence $\krn j_{\e K\!/\be
F}=C_{F}$, $\cok j^{\e\prime}_{\e K\!/\be F}=0$ and $\cok j_{\e
K\!/\be F}=C_{K}^{\e G}$ . The theorem is now immediate from
Theorem 2.4 and Proposition 2.2.\qed

\medpagebreak

We note that the preceding theorem applies, in particular, to ray
class fields over $F$ since such fields are equal to their own genus
field relative to $F\,$\footnote{ See [27, Lemma 2.1] for a general
description of the abelian extensions of $F$ which are equal to
their own genus field relative to $F$.}.

\medpagebreak

{\bf{Remark 4.2.}} There exists a function field analog of the
preceding theorem. Indeed, let $F$ be a function field, let $\infty$
be a fixed place of $F$ and let $S=\{\infty\}$. Further, let $K$ be
a finite abelian extension of $F$ where $\infty$ is tamely ramified
and the decomposition and inertia groups of $\infty$ agree. These
extensions were called ``of type (*)" in [7]. Assume further that
$K$ is its own $S$-genus field. Then the main theorem of [7],
combined with (5), shows that $(C_{K,S})^{\e G}_{\roman{trans}}=0$.
Consequently, Theorem 2.4 yields an exact sequence
$$
0\ra C_{F,S}\ra H^{1}\lbe(G,U_{K,S})\ra\bigoplus_{v\neq\infty}\Bbb
Z/e_{v}\e\Bbb Z\ra 0.
$$
In particular, if $h_{F,S}=1$, there exists an isomorphism
$$
H^{1}\lbe(G,U_{K,S})\simeq\bigoplus_{v\neq\infty} \Bbb
Z/e_{v}\e\Bbb Z.
$$
The above generalizes [7, Corollary 4.3].

\medpagebreak

\heading 5. Cyclic extensions.\endheading

In this Section we assume that $K/F$ is a cyclic extension of
degree $n$.

\medpagebreak

Set
$$
W_{F,S}=U_{F,S}\cap N_{K/F}K^{*}.\tag 22
$$
Then, by (4) and the periodicity of the cohomology of cyclic
groups, we have
$$\align
H^{1}(G,K^{*}\be/\e
U_{K,S})&\simeq\krn\ng\left[H^{2}(G,U_{K,S})\ra
H^{2}(G,K^{*})\right]\\
&\simeq\krn\ng\left[\e U_{F,S}\e/@,@,@, N_{\be K\be/\lbe F}\e
U_{\lbe K,S}\ra F^{\e *}\lbe/@,@,@, N_{\be K\be/\lbe F}\e
K^{*}\right]\\
&=W_{F,S}\e/\e N_{\be K\be/\lbe F}\e U_{\lbe K,S}
\endalign
$$
(cf. [22, Theorem 2(B), p.161]). Thus, by Proposition 2.2, the
following holds.

\proclaim{Lemma 5.1} Assume that $K/F$ is a cyclic extension of
global fields. Then there exists a natural exact sequence
$$
0\ra\cok j^{\e\prime}_{\e K\!/\be F,S}\ra\cok j_{\e K\!/\be
F,S}\ra W_{F,S}\e/\e N_{\be K\be/\lbe F}\e U_{\lbe K,S}\ra 0,
$$
where $W_{F,S}$ is the group (22).\qed
\endproclaim

\medpagebreak

The next result may be regarded as a ``dual" of Theorem 2.7.

\proclaim{Theorem 5.2} Let $K/F$ be a cyclic extension of global
fields. Assume that $U_{F,S}\subset N_{K/F}\e K^{*}$. Then $\cok
j_{\e K\!/\be F,S}$ contains at least
$$
\frac{\tprod_{\e v\in
S^{\e\prime}}d_{v}\big/D}{\left(n_{0},\tprod_{\e v\in
S^{\e\prime}}d_{v}/D\right)}
$$
classes, where $d_{v}$, $D$ and $n_{0}$ are the integers (9), (10)
and (11), respectively.
\endproclaim
{\it Proof.} The lemma and the proof of Theorem 2.7 immediately
yield the formula
$$
\frac{n}{\left(n,\tprod_{\e v\in S^{\e\prime}}d_{v}\right)} [\e
U_{F,S}\e\:\be W_{F,S}\e]\e\big[\cok j_{\e K\!/\be
F,S}\e\big]=\frac{\tprod_{\e v\in S^{\e\prime}}\e
d_{v}}{\left(n,\tprod_{\e v\in
S^{\e\prime}}d_{v}\right)}\left[\krn j_{\e K\!/\be F,S}\right].
$$
Therefore
$$
\frac{\tprod_{\e v\in S^{\e\prime}}\e d_{v}}{\left(n,\tprod_{\e
v\in S^{\e\prime}}d_{v}\right)}=\frac{\tprod_{\e v\in
S^{\e\prime}}d_{v}\big/D}{\left(n_{0},\tprod_{\e v\in
S^{\e\prime}}d_{v}/D\right)}
$$
divides $[\e U_{F,S}\e\:\be W_{F,S}\e]\e\big[\cok j_{\e K\!/\be
F,S}\e\big]$. The result is now clear since $W_{F,S}=U_{F,S}$ by
hypothesis.\qed

\medpagebreak

{\bf{Example 5.3.}} Let $F$ be either $\Bbb Q$, an imaginary
quadratic number field (with $S=S_{\infty}$ in both cases) or a
function field with $\#S=1$. Let $K$ be a cyclic extension of $F$
of degree $l^{\e m}$, where $m$ is a positive integer and $l$ is a
rational prime which is either $\geq 5$ in the number field case
or prime to $q-1$ in the function field case. Then $l^{\e m}$ is
prime to the order of the finite group $U_{F,S}$, whence
$U_{F,S}=N_{K/F}\e U_{K,S}\subset N_{K/F}\e K^{*}$. Thus the
hypothesis of Theorem 5.2 is satisfied. Assume now, for
simplicity, that in the function field case the prime in $S$
splits completely in $K$. Let $l^{t_{1}},l^{t_{2}},\dots,
l^{t_{r}}$ be the ramification indices of the ramified primes of
$K/F$, and assume that $m<\sum_{i=1}^{r}t_{i}$. Then the theorem
asserts that $l^{t}\!\mid\be[\e\cok j_{\e K\!/\be F,S}\e]$, where
$$
t=\sum_{i=1}^{r}t_{i}-m.
$$
In particular, if $K/k(t)$ is a real Artin-Schreir extension,
i.e., $K=F(\e\alpha\e)$, where $\alpha$ is a root of
$x^{p}-x+Q(t)=0$ and $Q(t)$ is such that the infinite prime of
$k(t)$ splits in $K$, then
$$
\roman{rank}_{\e p}\!\left(C_{K}^{\e G}\right)\geq r-1,
$$
where $r$ is the number of finite primes of $k(t)$ which ramify in
$K\,$\footnote{See Corollary 6.2 for the case of imaginary
Artin-Schreir extensions.}.

\medpagebreak

We conclude this Section by noting that (15) and the
identifications
$$
H^{1}(G,K^{*}\be/\e U_{K,S})\simeq W_{F,S}\e/\e N_{\be K\be/\lbe
F}\e U_{\lbe K,S}
$$
and
$$
H^{2}(G,U_{\be K,S}\be)\simeq U_{F,S}\e/\e N_{\be K\be/\lbe F}\e
U_{\lbe K,S}
$$
show that $U_{F,S}/\e W_{F,S}$ is isomorphic to a subgroup of
$\roman{B}@,@,@,({\Cal O}_{\be F,S},{\Cal O}_{\be K,S})$.
Consequently, Lemma 3.2 yields the following generalization of
[29, Theorem 1(iv)] (see also [op.cit., Lemma 5]).

\proclaim{Theorem 5.4} The index $[\e U_{F,S}\e\:\be W_{F,S}\e]$
divides $\left(\e\tprod_{\e v\in S^{\e\prime}}d_{v}\right)\lbe/D$.
Consequently, if the integers $d_{v}$, for $v\in S^{\e\prime}$,
are pairwise coprime, then every $S$-unit of $F$ is a norm from
$K$.
\endproclaim

\medpagebreak

{\bf Remark 5.5.} Note that in the ``limit case" $[\e
U_{F,S}\e\:\be W_{F,S}\e]=\left(\e\tprod_{\e v\in
S^{\e\prime}}d_{v}\right)\lbe/D$ (cf. Theorem 3.7), the proof of
Theorem 5.2 yields the identity
$$
\left[\krn j_{\e K\!/\be F,S}\right]=n_{@,@, 0}\be\left[\cok j_{\e
K\!/\be F,S}\right],
$$
where $n_{0}=n/D$. In particular, $\krn j_{\e K\!/\be F,S}$
contains at least $n_{@,@, 0}$ elements. Compare Remark 3.8.

\bigskip

\heading 6. Imaginary extensions of function fields.\endheading

The main result of this section is the following

\proclaim{Theorem 6.1} Let $F$ be a function field, let $K$ be a
Galois extension of $F$ of degree $n$ and let $k^{\e\prime}$ be
the field of constants of $K$. Assume that $\# S_{K}=1$ and that
$n$ is prime to $q^{\e\prime}-1$, where $q^{\e\prime}$ is the
cardinality of $k^{\e\prime}$. Then there exist an exact sequence
$$
0\ra C_{F,S}\ra C_{K,S}^{\e G}\ra\bigoplus_{R\setminus S}\Bbb
Z/e_{v}\e\Bbb Z\ra 0
$$
and an isomorphism
$$
H^{1}(G,C_{K,S})\simeq\krn\ng\left[\,\bigoplus_{v\in
S^{\e\prime}}\e\Bbb Z/d_{v}\e\Bbb
Z\overset{\Sigma}\to{\longrightarrow}\Bbb Z/ D\e\Bbb Z\e\right].
$$
\endproclaim
{\it Proof.} The hypothesis $\# S_{K}=1$ (i.e., $S$ consists of
exactly one prime of $F$ and there is only one prime of $K$ lying
above it) implies, by Dirichlet's Unit Theorem, that
$U_{K,S}=(k^{\e\prime})^{*}$. Hence
$H^{i}(G,U_{K,S})=H^{i}(G,\left(k^{\e\prime}\right)^{*})$ is
annihilated by multiplication by $(n,\e q^{\e\prime}\be-1)$ for
all $i\geq 1$. As $(n,\e q^{\e\prime}\be-1)=1$ by hypothesis, we
conclude that $H^{i}(G,U_{K,S})=0$ for all $i\geq 1$. The theorem
now follows by combining Theorems 2.4 and 3.3.\qed

\medpagebreak

The following result is one possible application of Theorem 6.1.

\proclaim{Corollary 6.2} Let $F=k(t)$ with $k$ a finite field of
characteristic $p$, let $Q(t)\in F^{*}$ and let $K=F(\e\alpha\e)$,
where $\alpha$ is a root of $x^{p}-x+Q(t)=0$. Assume that $K$ is
imaginary, i.e., that there exists only one prime of $K$ lying
above the infinite prime of $F$. Then there exists a canonical
isomorphism
$$
C_{K}^{\e G}\simeq(\Bbb Z/p\e\Bbb Z)^{r},
$$
where $r$ is the number of finite primes of $F$ which ramify in
$K$.
\endproclaim
{\it Proof.} Note that $K$ is a Galois extension of $F$ of degree
$p$, which is prime to $q^{\e\prime}-1$. Further, $h_{F}=1$. The
corollary is immediate from the theorem.\qed

\heading 7. Large $S$.\endheading

A (nonempty) set $S$ of primes of a global field $F$ is {\it large
relative to $K/F$} if $S$ contains all archimedean primes of $F$ and
all primes that ramify in $K/F$. In this section we assume that our
set $S$ is large. Note then that $S^{\e\prime}=S\cup(R\setminus
S)=S$.

\proclaim{Theorem 7.1} Let $K/F$ be a finite Galois extension of
global fields with Galois group $G$, and let $S$ be a set of
primes of $F$ which is large relative to $K/F$ (as defined above).
Then there exists an exact sequence
$$
0\ra H^{1}(G,U_{K,S})\ra C_{F,S}\ra C_{K,S}^{\e G}\ra{\cyr
X}^{2}(G,U_{K,S})\ra 0,
$$
where
$$
{\cyr X}^{2}(G,U_{K,S})=\krn\ng\left[\e
H^{2}(G,U_{K,S})\ra\bigoplus_{v\in S}
H^{2}(G_{w},\kwstar)\e\right],
$$
the map involved being induced by the inclusions
$U_{K,S}\subset\kwstar$ ($\, v\in S$).
\endproclaim
{\it Proof.} Since $S\supset R$, Theorem 2.4 shows that $\krn
j_{\e K\!/\be F,S}=H^{1}\lbe(G,U_{K,S})$ and that $\cok
j^{\e\prime}_{\e K\!/\be F,S}=0$. Now Proposition 2.2 and Remark
3.4 yield $\cok j_{\e K\!/\be F,S}={\cyr X}^{\e
2}\lbe(G,U_{K,S})$, where ${\cyr X}^{\e 2}(G,U_{K,S})$ is as in
the statement. \qed

\medpagebreak

The following corollary of the above theorem should be compared
with [9, Proposition A.2].

\proclaim{Corollary 7.2} Let $K/F$ be a finite Galois extension of
global fields with Galois group $G$. Assume that

\roster

\item"({\it a})" exactly one prime $v_{0}$ of $F$ ramifies in $K$,

\item"({\it b})" $S\supset S_{\infty}\cup\{v_{0}\}$ (i.e., $S$ is
large relative to $K/F$), and

\item"({\it c})" every prime in $S\setminus\{\e v_{0}\}$ splits
completely in $K$.

\endroster

Then there exists an exact sequence
$$
0\ra H^{1}(G,U_{K,S})\ra C_{F,S}\ra C_{K,S}^{\e G}\ra
H^{2}(G,U_{K,S})\ra 0.
$$
\endproclaim
{\it Proof.} By Lemma 3.2, $\roman{B}({\Cal O}_{\be F,S},{\Cal
O}_{\be K,S})$ is a group of order
$$
\left(\e\tprod_{v\in S}\e [K_{w}\be\:\be
F_{v}]\e\right)/\,\roman{l.c.m.}\!\left\{[K_{w}\be\:\be F_{v}]\:
v\in S\e\right\}=\left[ K_{w_{0}}\be\:\be
F_{v_{0}}\right]\lbe\big/\left[ K_{w_{0}}\be\:\be
F_{v_{0}}\right]=1.
$$
Consequently
$$
{\cyr X}^{\e 2}(G,U_{K,S})=\krn\ng\left[\e H^{2}(G,U_{K,S})\ra
\roman{B}({\Cal O}_{\be F,S},{\Cal O}_{\be K,S})\e\right]=
H^{2}(G,U_{K,S})
$$
(see Remark 3.4). The corollary is now immediate from the
theorem.\qed

\bigskip

\heading 8. Divisor class groups.\endheading

This Section may be regarded as the ``$S=\emptyset$" version of
Sections 2 and 3 in the function field case. We follow [22]
closely.

\medpagebreak

Let $K/F$ be a Galois extension of degree $n$ of function fields
with Galois group $G$. Let $X^{\prime}$ denote the unique smooth
complete curve over $k^{\e\prime}$ with function field $K$.
Similarly, let $X$ be the smooth complete curve over $k$ with
function field $F$. Further, let $H=\gal(K/Fk^{\e\prime})$ and
$g=\gal(k^{\e\prime}\lbe/k)$, so that $G/H=g$. The group of
$k^{\e\prime}$-rational points of the Jacobian variety of
$X^{\e\prime}$ will be denoted by $J_{K}$. Since $X^{\prime}$ has a
$k^{\e\prime}$-rational point, we have
$J_{K}=\roman{Pic}^{0}(X^{\e\prime})$. Similarly, let
$J_{F}=\roman{Pic}^{0}(X)$. Then there exists a natural {\it
capitulation map} $j_{\e K\be/\be F}\:J_{F}\ra J_{K}^{\e G}$, which
is defined by pulling back line bundles. Its kernel was determined
in [22, Theorem 5] under the assumption that $k=k^{\e\prime}$. In
general, the following holds

\proclaim{Theorem 8.1} Let $K/F$ be a finite Galois extension of
function fields and let $M$ be the maximal abelian unramified
extension of $k^{\e\prime}\lbe F$ in $K$. Then the kernel of the
capitulation map $j_{\e K\be/\be F}\:J_{F}\ra J_{K}^{\e G}$ is
naturally isomorphic to $\roman{Hom}_{g}(\roman{Gal}(\lbe
M/Fk^{\e\prime}\e), (k^{\e\prime})^{*})$.
\endproclaim
{\it Proof.} The Hochschild-Serre spectral sequence
$$
H^{p}(\e g,H^{q}(H,(k^{\e\prime})^{*}))\implies
H^{p+q}(G,(k^{\e\prime})^{*})
$$
yields an exact sequence
$$
0\ra H^{1}(g,(k^{\e\prime})^{*})\ra H^{1}(G,(k^{\e\prime})^{*})\ra
H^{1}(H,(k^{\e\prime})^{*})^{\e g}\ra\br(k^{\e\prime}\lbe/\e k\e)
$$
(see, for example, [21, p.309]). Now
$H^{1}(g,(k^{\e\prime})^{*})=0$ by Hilbert's Theorem 90 and
$\br(k^{\e\prime}\lbe/\e k\e)=0$ since the Brauer group of a
finite field is zero (see, e.g., [24]). Therefore
$$\align
H^{1}(G,(k^{\e\prime})^{*})&\simeq H^{1}(H,(k^{\e\prime})^{*})^{\e
g}\simeq\roman{Hom}_{g}(H,(k^{\e\prime})^{*})\\
&\simeq\roman{Hom}_{g}(\roman{Gal}(\lbe
M^{\e\prime}/Fk^{\e\prime}\e), (k^{\e\prime})^{*}),
\endalign
$$
where $M^{\e\prime}$ is the maximal abelian extension of
$Fk^{\e\prime}$ in $K$. Taking into account these facts, it is not
difficult to adapt the proof of Theorem 5 of [22] to the case
where $k^{\e\prime}$ is not necessarily equal to $k$. We leave the
details to the reader.\qed

\medpagebreak

We will now study the cokernel of the capitulation map $J_{F}\ra
J_{K}^{\e G}$. Let ${\s I}_{K}^{\e
0}=\roman{Div}^{0}(X^{\e\prime})$. Then there exists a natural
exact sequence
$$
1\ra (k^{\e\prime})^{*}\ra K^{*}\ra{\s I}_{K}^{\e 0}\ra J_{K}\ra
0\tag 23
$$
inducing an exact sequence
$$
1\ra (K^{*}/(k^{\e\prime})^{*})^{G}\ra\left({\s I}_{K}^{\e
0}\right)^{\lbe G}\ra J_{K}^{\e G}\ra
H^{1}(G,K^{*}/(k^{\e\prime})^{*})\ra H^{1}\ng\left(G,{\s
I}_{K}^{\e 0}\right).\tag 24
$$
Define
$$
\left(J_{K}^{\e G}\right)_{\roman{trans}}=\krn\ng\left[\e
J_{K}^{\e G}\ra H^{1}(G,K^{*}/(k^{\e\prime})^{*})\right].
$$
Then the image of $j_{\e K\be/\be F}$ is contained in
$\left(J_{K}^{\e G}\right)_{\roman{trans}}$. Let $j^{\e\prime}_{\e
K\be/\be F}\:J_{F}\ra \left(J_{K}^{\e G}\right)_{\roman{trans}}$
be the map induced by $j_{\e K\be/\be F}$. Arguing as in \S 2
(see, especially, diagram (8)), we obtain the following result.

\proclaim{Proposition 8.2} There exists a natural exact sequence
$$
0\ra\krn j_{\e K\be/\be F}\ra H^{1}(G,(k^{\e\prime})^{*})\ra
\left({\s I}_{K}^{\e 0}\right)^{\lbe G}\ng\big/\e{\s I}_{F}^{\e
0}\ra\cok j^{\e\prime}_{\e K\be/\be F}\ra 0.\qed
$$
\endproclaim

\medpagebreak

Now define positive integers $\delta$ and
$\delta^{\e\prime}\,$\footnote{One might call these integers the
``index" and the ``period" of $K/F$, respectively.} by
$$\align
\img\ng\be\left[\e\roman{Div}(X^{\e\prime})^{G}\overset{\roman{deg}}\to\longrightarrow\Bbb
Z\e\right]&=\delta\e\Bbb
Z\\
\img\ng\be\left[\e\roman{Pic}(X^{\e\prime})^{G}\overset{\roman{deg}}\to\longrightarrow\Bbb
Z\e\right]&=\delta^{\e\prime}\e\Bbb Z.
\endalign
$$
The integer $\delta$ is ``fairly easy to compute" [22, p.164].
Indeed, it can be shown that, if $P_{1},P_{2},\dots,P_{r}$ are the
primes of $F$ that ramify in $K$ and $e_{1},e_{2},\dots, e_{r}$
are their respective ramification indices, then
$$
\delta=(n,(n/e_{1})\e\roman{deg}_{F}P_{1},\dots,
(n/e_{r})\e\roman{deg}_{F}P_{r}).
$$
See [22, proof of Proposition 1, p.163] (by contrast,
$\delta^{\e\prime}$ is a ``more subtle invariant" [op.cit.,
p.164]). Clearly, $\delta^{\e\prime}\!\mid\!\delta\!\mid\! n$.

\proclaim{Lemma 8.3} There exists a natural exact sequence
$$
0\ra\left({\s I}_{K}^{\e 0}\right)^{\lbe G}\ng\big/\e{\s
I}_{F}^{\e 0}\ra\bigoplus_{i=1}^{r}\Bbb Z/e_{i}\e\Bbb
Z\ra\delta\,\Bbb Z/n\,\Bbb Z\ra 0.
$$
\endproclaim
{\it Proof.} This follows at once from [22, Proposition 1 and
Theorem 3(B)].\qed

\medpagebreak

We also have the following exact sequence, which is the analog of
Proposition 2.2 in this context.

\proclaim{Proposition 8.4} There exists an exact sequence
$$
0\ra\cok j^{\e\prime}_{\e K\be/\be F}\ra\cok j_{\e K\be/\be F}\ra
H^{1}(G,K^{*}/(k^{\e\prime})^{*})\ra\g\ra 0,
$$
where $\g$ is a cyclic group of order $\delta/\delta^{\e\prime}$.
\endproclaim
{\it Proof.} This is immediate from (24) and [22, Theorem 3(A),
p.163], which shows that the image, $\g$, of
$H^{1}(G,K^{*}/(k^{\e\prime})^{*})$ in $H^{1}\!\left(G,{\s
I}_{K}^{\e 0}\right)$ under the map appearing in (24) is a cyclic
group of the indicated order. \qed

\medpagebreak

Combining statements 8.2, 8.3 and 8.4, we obtain the following
analogue of C.Chevalley's classical ``ambiguous class number
formula" (cf. [22, Theorem 8, p.166]):

\proclaim{Theorem 8.5} We have
$$
\frac{\left[\e J_{K}^{G}\e\right]}{\left[\e J_{F}\e\right]}=
\frac{\left[\e H^{1}(G,K^{*}\be/\e
(k^{\e\prime})^{*})\e\right]\prod_{\e i=1}^{\e r}\lbe
e_{i}}{(n/\delta^{\e\prime})\left[\e
H^{1}(G,(k^{\e\prime})^{*})\e\right]}.\qed
$$
\endproclaim

\medpagebreak

Now we note that (23) induces an exact sequence
$$
0\ra H^{1}(G,K^{*}\be/\e (k^{\e\prime})^{*})\ra
H^{2}(G,(k^{\e\prime})^{*})\ra\br(K/F)\ra H^{2}(G,K^{*}\be/\e
(k^{\e\prime})^{*})
$$
(cf. (4)). On the other hand, there exists a natural map
$\br(K/F)\ra H^{2}(G,{\Cal I}_{K})$ (where ${\Cal
I}_{K}=\iv(X^{\e\prime})$) which factors as
$$
\br(K/F)\ra H^{2}(G,K^{*}\be/\e (k^{\e\prime})^{*})\ra
H^{2}(G,{\Cal I}_{K}),
$$
where the second map is induced by the canonical map $K^{*}\be/\e
(k^{\e\prime})^{*}\ra {\Cal I}_{K}$. We conclude that there exists
an exact sequence
$$
0\ra H^{1}(G,K^{*}\be/\e (k^{\e\prime})^{*})\ra
H^{2}(G,(k^{\e\prime})^{*})\ra B,\tag 25
$$
where
$$
\roman{B}=\krn\ng\left[\,\br(K/F)\ra H^{2}(G,{\Cal I}_{K})\e\right].
$$
Now essentially the same argument which proves Lemma 3.2 yields an
isomorphism
$$
B\simeq\krn\ng\left[\e\bigoplus_{i=1}^{r}\e\Bbb Z/\e d_{i}\e\Bbb
Z\overset{\Sigma}\to{\longrightarrow}\Bbb Z/\e[\e
d_{1},d_{2},\dots,d_{r}]\e\Bbb Z\e\right],\tag 26
$$
where $d_{i}=[\e H^{2}(G_{P_{i}},U_{P_{i}})\e]$ for each
$i=1,2,\dots,r$ and $[\e d_{1},d_{2},\dots,d_{r}]$ denotes the
least common multiple of $d_{1},d_{2},\dots,d_{r}$. In particular,
$B$ is a group of order $\left(\e\prod_{i=1}^{r}
d_{i}\e\right)\be/[\e d_{1},d_{2},\dots,d_{r}]$. Consequently, if
$K/F$ is cyclic, we have
$$
\left[B\e\right]=\frac{\prod_{i=1}^{r} e_{i}}{[\e
e_{1},e_{2},\dots,e_{r}]}\,.
$$
Now define
$$
m(P_{i})=\left(\frac{q^{\e\roman{deg}_{F}
P_{i}}-1}{\left(q^{\e\roman{deg}_{F}
P_{i}}-1,e_{i}\right)},q-1\right)\qquad(1\leq i\leq r)
$$
and set
$$
m=\left(m(P_{1}),m(P_{2}),\dots,m(P_{r})\right).\tag 27
$$
The following proposition shows that the degrees and the
ramification indices of the primes that ramify in a cyclic extension
$K/F$ are subject to certain non-obvious constraints.

\proclaim{Proposition 8.6} Assume that $K/F$ is a cyclic
extension. Then
$$
\frac{\prod_{i=1}^{r} e_{i}}{[\e e_{1},e_{2},\dots,e_{r}]}\,\equiv
0\,\left(\be\roman{mod}\,\frac{q-1}{m}\right),
$$
where $m$ is the integer (27). Consequently, if the ramification
indices $e_{i}$ are pairwise coprime, then
$$
\frac{q^{\e\roman{deg}_{F} P_{i}}-1}{\left(q^{\e\roman{deg}_{F}
P_{i}}-1,e_{i}\right)}\equiv 0\ng\pmod{q-1}
$$
for every $i=1,2,\dots,r$.
\endproclaim
{\it Proof.} The periodicity of the cohomology of cyclic groups
and [22, Theorem 2(B), p.161, and Proposition 2, p.166] show that
$$
\frac{\left[\e H^{2}(G,(k^{\e\prime})^{*})\e\right]}{\left[\e
H^{1}(G,K^{*}\be/\e (k^{\e\prime})^{*})\e\right]}=\frac{q-1}{m}\,.
$$
Therefore, $(q-1)/m$ divides $[\e B\e]=\left(\e\prod_{i=1}^{r}
e_{i}\right)/\e [\e e_{1},e_{2},\dots,e_{r}]$, as claimed. Now, if
the ramification indices $e_{i}$ are pairwise coprime, then
necessarily $m=q-1$. The definition of $m$ now yields the second
assertion of the proposition.\qed

\medpagebreak

We conclude this paper with the following strengthening of [22,
Theorem 14].

\proclaim{Theorem 8.7} Assume that $K/F$ is a Galois extension
with Galois group $G\simeq\bigoplus_{i=1}^{s}\be\Bbb Z/\ell\e\Bbb
Z$, where $s$ is a positive integer and $\ell$ is a prime which
divides $q-1$. Assume, furthermore, that the field of constants of
$F$ is algebraically closed in $K$. Then
$$
\roman{rank}_{\e\ell}\!\left(J_{K}^{\e G}\right)\geq
(s(s+1)/2)-r\,.
$$
\endproclaim
{\it Proof.} By (26), $\roman{rank}_{\e\ell}(B)\leq r-1$.
Consequently, (25) shows that
$$
\roman{rank}_{\e\ell}(H^{2}(G,(k^{\e\prime})^{*}))=
\roman{rank}_{\e\ell}(H^{2}(G,k^{*}))\leq
\roman{rank}_{\e\ell}(H^{1}(G,K^{*}\be/\e k^{*}))+r-1.
$$
The rest of the proof is similar to that of Theorem 14 of
[22].\qed

\medpagebreak

\heading Appendix.\endheading

The integers $\left[H^{2}\be(G_{w},U_{\be w})\right]$, where $w$
lies above a ramified prime of $K/F$, intervene at various places
in the paper. The following result (which may be well-known)
relates these integers to the ramification indices $e_{v}$ of
$K/F$.

\proclaim{Proposition A.1} Let $v$ be a non-archimedean prime of
$F$ and let $w$ be a fixed prime of $K$ lying above $v$. Then
there exists an exact sequence
$$
0\ra \Bbb Z/(e_{v},f_{v})\e\Bbb Z\ra H^{2}(G_{w},U_{w})\ra \Bbb
Z/e_{v}\e\Bbb Z,
$$
where $f_{v}$ is the residue class degree. In particular,
$\left[H^{2}\be(G_{w},U_{w})\right]$ divides $ e_{v}^{2}$.
\endproclaim
{\it Proof}. Let $F_{v}^{\e\roman{un}}$ be the maximal unramified
extension of $F_{v}$ contained in $K_{w}$. Set
$I_{w}=\gal(K_{w}/F_{v}^{\e\roman{un}})$. Then $I_{w}$ is a
subgroup of $G_{w}$ of order $e_{v}$ and $G_{w}/\e I_{w}=\gal(
F_{v}^{\e\roman{un}}/F_{v})$ is a cyclic group of order $f_{v}$.
The exact sequence of terms of low degree belonging to the
Hochschild-Serre spectral sequence $H^{p}(G_{w}/\e
I_{w},H^{q}(I_{w},U_{w}))\implies H^{p+q}(G_{w},U_{w})$ yields,
since $H^{i}(G_{w}/\e I_{w},U_{w}^{I_{w}})=0$ for all $i\geq 1$,
an exact sequence
$$
0\ra H^{1}(G_{w}/\e I_{w},H^{1}(I_{w},U_{w}))\ra
H^{2}(G_{w},U_{w})\ra H^{2}(I_{w},U_{w})
$$
(see [21, p.309]). On the other hand, since
$K_{w}/F_{v}^{\e\roman{un}}$ is totally ramified, we have
$H^{2}(I_{w},U_{w})\simeq H^{2}(I_{w},K_{w}^{*})\simeq \Bbb
Z/e_{v}\e\Bbb Z$, by [24, Ch.XII, Exer.2(b), p.182] and local
class field theory. Further, $H^{1}(I_{w},U_{w})=\Bbb
Z/e_{v}\e\Bbb Z$ with trivial $G_{w}\be/ I_{w}$-action [op.cit.],
whence
$$
H^{1}(G_{w}/\e I_{w},H^{1}(I_{w},U_{w}))\simeq\roman{Hom}(\Bbb
Z/f_{v}\e\Bbb Z,\Bbb Z/e_{v}\e\Bbb Z)\simeq \Bbb
Z/(e_{v},f_{v})\e\Bbb Z.
$$
We conclude that there exists an exact sequence
$$
0\ra \Bbb Z/(e_{v},f_{v})\e\Bbb Z\ra H^{2}(G_{w},U_{w})\ra \Bbb
Z/e_{v}\e\Bbb Z,
$$
as asserted.\qed

\medpagebreak

Now define integers $d_{v}^{\,\prime}$ by
$$
d_{v}^{\,\prime}=\cases[K_{w}\:\be F_{v}]\quad\text{ if
$v\in S$}\\
\,\,\,\,\,\,\,\,e_{v}\,\,\,\qquad\e\text{ if $v\in R\setminus S$}.
\endcases
$$
The following corollary shows that the hypothesis of Corollary 3.5
is satisfied if the integers $d_{v}^{\,\prime}$ defined above (which
are relatively easy to compute) are pairwise coprime.

\proclaim{Corollary A.2} Assume that the integers
$d_{v}^{\,\prime}$, where $v\in S^{\e\prime}$, are pairwise
coprime. Then so also are the integers $d_{v}$ defined by
(18).\qed

\endproclaim

\enddocument